\documentclass[graybox]{svmult}

\usepackage{type1cm}        
\usepackage{makeidx}         
\usepackage{graphicx}        
\usepackage{multicol}        
\usepackage[bottom]{footmisc}
\usepackage[numbers]{natbib}
\usepackage{orcidlink}

\usepackage{newtxtext}       %
\usepackage[varvw]{newtxmath}       

\usepackage{amsmath,graphicx,textcomp,hyperref,tabularx,xcolor,theorem}

\hypersetup{
	colorlinks=true, 
	linkcolor=blue, 
	urlcolor=red, 
	citecolor=magenta,
	linktoc=all 
}

\usepackage{algorithmicx}
\usepackage{algorithm}
\usepackage{algpseudocode}

\usepackage{xpatch}
\makeatletter
\xpatchcmd{\algorithmic}
{\ALG@tlm\z@}{\leftmargin\z@\ALG@tlm\z@}
{}{}
\makeatother

\newenvironment{tmparmod}[3]{\begin{list}{}{\setlength{\topsep}{0pt}\setlength{\leftmargin}{#1}\setlength{\rightmargin}{#2}\setlength{\parindent}{#3}\setlength{\listparindent}{\parindent}\setlength{\itemindent}{\parindent}\setlength{\parsep}{\parskip}} \item[]}{\end{list}}

\newcounter{tmcounter}

\usepackage[T3,T1]{fontenc}
\DeclareSymbolFont{tipa}{T3}{cmr}{m}{n}
\DeclareMathAccent{\invbreve}{\mathalpha}{tipa}{16}

\newcommand{\myvector}[1]{\text{{\hspace{0.1em}}#1{\hspace{0.15em}}}}

\newcommand{\vD}{\myvector{D}}

\newcommand{\ba}{\boldsymbol{a}}
\newcommand{\bI}{\boldsymbol{I}}
\newcommand{\bb}{\boldsymbol{b}}

\newcommand{\bc}{\boldsymbol{c}}

\newcommand{\bzero}{\boldsymbol{0}}
\newcommand{\bS}{\boldsymbol{S}}

\newcommand{\fA}{\boldsymbol{A}}
\newcommand{\fB}{\boldsymbol{A}}

\newcommand{\emh}{{e - \frac{1}{2}}}

\newcommand{\eph}{{e + \frac{1}{2}}}

\newcommand{\uu}{\boldsymbol{u}}
\newcommand{\uud}{\uu^{\delta}}

\newcommand{\uudhat}{\hat{\uu}^{\delta}}
\newcommand{\discf}{\ensuremath{\pf_h^\delta}}
\newcommand{\discfref}{\ensuremath{\hat{\pf}_{h,e}^\delta}}
\newcommand{\discfrefe}{\ensuremath{\hat{\pf}_h^\delta}}
\newcommand{\discF}{\ensuremath{\F_h^\delta}}
\newcommand{\discFref}{\ensuremath{\hat{\F}_{h,e}^\delta}}
\newcommand{\discFrefe}{\ensuremath{\hat{\F}_h^\delta}}
\newcommand{\dfrx}{\partial_x^\text{FR}}
\newcommand{\dlocx}{\partial_x^\text{loc}}
\newcommand{\pdx}{\partial_x}
\newcommand{\fnum}{\pf^\text{num}}

\newcommand{\en}{\ensuremath{\mathbb{E}}}

\newcommand{\uU}{\boldsymbol{U}}

\newcommand{\pf}{\boldsymbol{f}}
\newcommand{\fmatrix}{\ensuremath{\mathcal{A}}}
\newcommand{\fp}{\pf^+}
\newcommand{\fm}{\pf^-}

\newcommand{\F}{\boldsymbol{F}}

\newcommand{\bss}{\boldsymbol{s}}

\newcommand{\paragraphtoc}[1]{}
\providecommand{\citep}{}
\providecommand{\citet}{}

\newcommand{\re}{\mathbb{R}}

\newcommand{\poly}{\mathbb{P}}

\newcommand{\ud}{\textrm{d}}

\newcommand{\half}{\frac{1}{2}}

\newcommand{\Fnum}{\F^\text{num}}

\newcommand{\nc}{M}

\newcommand{\bv}{\boldsymbol{v}}

\usepackage{url}

\makeindex             

\begin{document}

\title*{Jin-Xin relaxation as a shock-capturing method for high-order DG/FR schemes}
\titlerunning{Jin-Xin relaxation as a shock-capturing method}
\author{
Marco Artiano\orcidID{0000-0003-1903-4107}, \\
Arpit Babbar\orcidID{0000-0002-9453-370X},\\
Michael Schlottke-Lakemper \orcidID{0000-0002-3195-2536}, \\
Gregor Gassner \orcidID{0000-0002-1752-1158}, and \\
Hendrik Ranocha\orcidID{0000-0002-3456-2277}
}
\authorrunning{Artiano, Babbar, Schlottke-Lakemper, Gassner, Ranocha}
\institute{
Marco Artiano \at Institute of Mathematics, Johannes Gutenberg University, Mainz, \email{martiano@uni-mainz.de}
\and Arpit Babbar \at Institute of Mathematics, Johannes Gutenberg University, Mainz, \email{ababbar@uni-mainz.de}
\and Michael Schlottke-Lakemper \at High-Performance Scientific Computing, Centre for Advanced Analytics and Predictive Sciences, University of Augsburg, Germany, \email{michael.schlottke-lakemper@uni-a.de},
\and Gregor Gassner \at Mathematical Institute, University of Cologne, Cologne, Germany \email{ggassner@math.uni-koeln.de}
\and Hendrik Ranocha \at Institute of Mathematics, Johannes Gutenberg University, Mainz, \email{hendrik.ranocha@uni-mainz.de}
}
\maketitle
\abstract*{
Jin-Xin relaxation is a method for approximating non-linear hyperbolic conservation laws by a linear system of hyperbolic equations with an $\varepsilon$ dependent stiff source term.
The system formally relaxes to the original conservation law as $\varepsilon \to 0$.
An asymptotic analysis of the Jin-Xin relaxation system shows that it can be seen as a convection-diffusion equation with a diffusion coefficient that depends on the relaxation parameter $\varepsilon$.
This work makes use of this property to use the Jin-Xin relaxation system as a shock-capturing method for high-order discontinuous Galerkin (DG) or flux reconstruction (FR) schemes.
The idea is to use a smoothness indicator to choose the $\varepsilon$ value in each cell, so that we can use larger $\varepsilon$ values in non-smooth regions to add extra numerical dissipation.
We show how this can be done by using a single stage method by using the compact Runge-Kutta FR method that handles the stiff source term by using IMplicit-EXplicit Runge-Kutta (IMEX-RK) schemes.
Numerical results involving Burgers' equation and the compressible Euler equations are shown to demonstrate the effectiveness of the proposed method.}

\section{Introduction}
Shock-capturing methods are important for high-order discretizations of hyperbolic conservation laws with non-smooth solutions.
The general idea is to add some form of numerical dissipation to prevent spurious oscillations.
For high-order discontinuous Galerkin (DG) methods, it has been studied since the initial works of Cockburn and Shu~\cite{Cockburn1989a}.
Their first approach was to use a TVB (total variation bounded) limiter~\cite{Cockburn1989a}, which reduces the DG polynomial to a constant or a linear function in \textit{troubled cells}.
Although the TVB limiter succeeds in preventing spurious oscillations, it is known to be too dissipative as it loses most of the information of the high-order scheme.
There have thus been several works that add numerical dissipation more carefully, such as moment limiters~\cite{Biswas1994}, WENO limiters~\cite{Qiu2005}, artificial viscosity~\cite{Persson2006}, and subcell-based limiters~\cite{hennemann2021}.
A recent review can be found in the introduction of~\cite{babbar2024admissibility}.

We propose to use the Jin-Xin relaxation system~\cite{jinxin1995} as a shock-capturing method for high-order DG/flux reconstruction (FR) schemes.
Based on the \emph{Chapman-Enskog} expansion of the Jin-Xin relaxation system~\cite{liu1987}, it can be interpreted as a convection-diffusion equation with a diffusion coefficient that depends on the relaxation parameter $\varepsilon$.
The idea is to use the smoothness indicator of~\cite{hennemann2021}, motivated by~\cite{Persson2006}, to choose the value of $\varepsilon$ in each cell, so that we add more numerical dissipation in troubled cells.
We show how this can be done in a single-stage time integration method by using the compact Runge-Kutta flux reconstruction (RKFR) method of~\cite{babbar2025crk}.
The approach is compatible with standard Runge-Kutta (RK) methods as well.
The numerical results are shown comparing the scheme with the subcell based blending limiter of~\cite{hennemann2021}, and the results show comparable resolution of the two methods.
At this stage, the issue with the proposed method is that it requires fine-tuning of the parameter $\varepsilon_{\max}$ that controls the maximum amount of numerical dissipation added in troubled cells, resolving which will be an interesting direction for future research.
The advantage of the method is that it does not require a non-linear Riemann solver.

In Section~\ref{sec:crkfr}, we briefly review the compact RKFR method for the time discretization of a hyperbolic conservation law with linear fluxes and stiff source terms.
This kind of system comes out of the Jin-Xin relaxation method discussed in Section~\ref{sec:jin.xin}, where we also show how it can be seen as a convection-diffusion equation with a diffusion coefficient that depends on the relaxation parameter $\varepsilon$.
In Section~\ref{sec:shock.capturing}, we describe how to use the Jin-Xin relaxation system as a shock-capturing method by using a smoothness indicator to choose the value of $\varepsilon$ in each cell.
Section~\ref{sec:num.results} shows some numerical results for scalar equations and the compressible Euler equations combined with the Jin-Xin relaxation approach as a shock-capturing method.
Finally, Section~\ref{sec:conclusion} concludes the paper with a summary and a future outlook.

\section{Compact Runge-Kutta FR method for stiff source terms} \label{sec:crkfr}

The Jin-Xin relaxation system described in detail in Section~\ref{sec:jin.xin} below is a hyperbolic system with a linear flux and stiff source terms. Thus, it requires a time discretization that can handle stiff source terms.
Such a system is generally written for 1-D as
\begin{equation}
\label{eq:lin.con.law.source} \uu_t + \pf(\uu)_x = \bss(\uu), \qquad \pf(\uu) = \fmatrix \uu, \quad \uu (x, t_0) = \uu_0
  (x), \quad x \in \Omega,
\end{equation}
where the solution $\uu(x,t) \in \re^{\nc}$ is the vector of conserved quantities, $\Omega$ is the physical domain, $\fmatrix$ is a matrix operator, $\bss(\uu)$ is the stiff source term demanding smaller time steps when used with explicit methods (see, e.g., Appendix 1 of~\cite{babbar2025crknoncons}), $\uu_0$ is the initial condition, and some boundary conditions are additionally prescribed.
In this section, we review the IMplicit-EXplicit (IMEX) compact Runge-Kutta flux reconstruction (cRKFR) scheme of~\cite{babbar2025crknoncons} for~\eqref{eq:lin.con.law.source}, although the scheme was originally proposed for general non-linear fluxes.
The IMEX cRKFR method treats the flux term explicitly and the source term implicitly.
The implicit equation for the source term is local to each solution point, and is thus efficient as it, e.g., does not require solving a global system of equations.

\begin{figure}
\begin{center}
\begin{tabular}{cc}
\includegraphics[width=0.45\columnwidth]{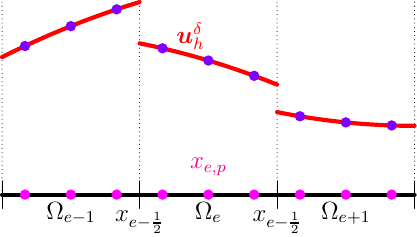} &
\includegraphics[width=0.45\columnwidth]{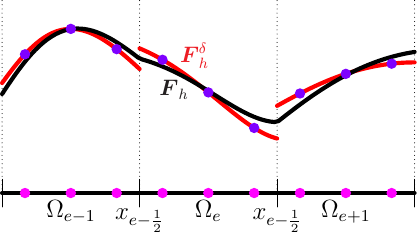}\\
(a) & (b)
\end{tabular}
\end{center}
\caption{(a) Piece-wise polynomial solution at time $t_n$, and (b)
discontinuous and continuous flux. The figure is inspired from~{\cite{babbar2022lax}}.}
\label{fig:solflux1}
\end{figure}

We begin by reviewing the notations for the finite element grid following~\cite{babbar2022lax,babbar2025crk}.
The physical domain $\Omega$ is divided into disjoint elements $\{\Omega_e \}$ with
\begin{equation}
  \label{eq:reference.element} \Omega_e = [x_{\emh}, x_\eph] \qquad
  \text{and} \qquad \Delta x_e = x_\eph - x_{\emh} .
\end{equation}
Each element is mapped to a reference element, $\Omega_e \to [0, 1]$, by
\begin{equation}
x \mapsto \xi = \frac{x - x_{\emh}}{\Delta x_e} .
\label{eq:ref.map}
\end{equation}
Inside each element, we approximate the numerical solution
to~\eqref{eq:lin.con.law.source} by $\poly_N$ functions, which are
polynomials of degree $N \geq 0$, so that the numerical solution is in the space
\begin{equation}
  \label{eq:fr.basis} V_h = \{v_h : v_h |_{\Omega_e} \in \poly_N \},
\end{equation}
which allows the solution to be discontinuous at element interfaces, see Figure~\ref{fig:solflux1}(a).

For each $v_h \in V_h$~\eqref{eq:fr.basis} and a physical element $\Omega_e$~\eqref{eq:reference.element}, we define its representative in the reference element $\hat{v}_{h, e} = \hat{v}_{h, e} (\xi)$ through the reference map~\eqref{eq:ref.map} as
\begin{equation}
  \label{eq:vh.reference} \hat{v}_{h, e} (\xi) = v_h  (x_{\emh} + \xi \Delta
  x_e) = v_h (x) .
\end{equation}
The element index $e$ will often be suppressed for brevity.
We represent $\poly_N$ polynomials by their values at the \textit{solution points} $0 \le \xi_0 < \xi_1 < \cdots < \xi_N \le 1$, which will be Gauss-Legendre (GL) or Gauss-Lobatto-Legendre (GLL) nodes in $[0, 1]$.
There are associated quadrature weights $w_j$ such that the quadrature rule is exact for
polynomials of degree up to $2 N + 1$ for GL points and up to degree $2 N - 1$
for GLL points.
Note that the nodes and weights we use are with respect to the interval $[0, 1]$ whereas they are usually defined for the interval $[-1,+1]$.
The numerical solution $\uud_h \in V_h^{\nc}$, where $\nc$ is the number of conservative variables~\eqref{eq:lin.con.law.source}, inside an element $\Omega_e$ is given in reference coordinates as
\begin{equation}
  \label{eq:soln.poly} x \in \Omega_e : \qquad \uudhat_{h, e} (\xi, t) =
  \sum_{p = 0}^N \uu_{e, p} (t) \ell_p (\xi),
\end{equation}
where each $\ell_p$ is a Lagrange polynomial of degree $N$ given by
\begin{equation}
  \label{eq:defn.lagrange} \ell_p (\xi) = \prod_{q = 0, q \ne p}^N \frac{\xi -
  \xi_q}{\xi_p - \xi_q} \in \poly_N, \quad \text{i.e.,} \quad \ell_p (\xi_q) = \delta_{p q} .
\end{equation}
The spatial derivatives required to solve~\eqref{eq:lin.con.law.source} are computed on the reference
interval $[0, 1]$ using the differentiation matrix $\vD = [D_{ij}]$ whose
entries are given by
\begin{equation}
  \label{eq:diff.matrix} D_{ij} = \ell_j' (\xi_i), \qquad 0 \le i, j \le N.
\end{equation}
Figure~\ref{fig:solflux1}(a) illustrates a piece-wise polynomial solution at some
time $t_n$ with discontinuities at the element boundaries.

A crucial ingredient of the cRKFR scheme and the FR scheme in general is the degree $N$ \textit{discontinuous flux} approximation $\discf$ (Figure~\ref{fig:solflux1}(b)), which is defined in reference coordinates for element $e$ as
\begin{equation} \label{eq:discts.flux}
\discfref(\xi) = \sum_{p=0}^N \pf(\uu_{e,p}) \ell_p(\xi).
\end{equation}
Then, the FR differentiation operator is defined as
\begin{equation}
\begin{gathered}
\dfrx \pf (\uud_h) = \pdx \pf_h, \quad \pf_h = \discf +
(\fnum_\eph - \fm_\eph) g_R + (\fnum_{\emh} - \fp_\emh) g_L, \\
\fm_\eph = \discfrefe (1), \quad \fp_\emh = \discfrefe (0),
\label{eq:dfrx.defn}
\end{gathered}
\end{equation}
where $g_L, g_R \in \poly_{N + 1}$ are FR correction functions~\cite{Huynh2007,Vincent2011a} satisfying $g_L (0) = g_R (1) = 1$, $g_L (1) = g_R (0) = 0$, and $\fnum_\eph$ is the numerical flux at the interface $x_\eph$. We use the Rusanov's flux~\cite{rusanov1962}, which for the linear flux $\pf(\uu) = \fmatrix \uu$ is given by
\begin{equation} \label{eq:num.flux.fr}
\fnum_\eph(\uud_h) = \half(\pf(\uu_\eph^-) + \pf(\uu_\eph^+)) - \frac {\sigma(\fmatrix)}{2} (\uu_\eph^+ - \uu_\eph^-), \qquad \uu_\eph^\pm = \uud_h(x_\eph^\pm),
\end{equation}
where $\sigma(\fmatrix)$ is the spectral radius of $\fmatrix$.
The $\pf_h$ defined in~\eqref{eq:dfrx.defn} is called the \textit{continuous
flux approximation} in the FR literature~{\cite{Huynh2007}} as it is globally
continuous, taking the numerical flux value $\fnum_{e + 1 / 2}$ at each
element interface $e + 1 / 2$.
We choose correction functions $g_{L/R} \in
\mathbb{P}_{N + 1}$ to be $g_{\text{Radau}}$ for GL and $g_2$ for GLL nodes~{\cite{Huynh2007}}.
A semi-discretization of~\eqref{eq:lin.con.law.source} is thus obtained as
\begin{equation}
  \label{eq:semidiscretization.dg} \partial_t \uud_h + \dfrx \pf (\uud_h) = \bss(\uud_h).
\end{equation}
The extension to multiple space dimensions is performed by using tensor products.
Solving~\eqref{eq:semidiscretization.dg} using an RK method gives the original RKFR method of~{\cite{Huynh2007}}.
We now describe the fully discrete cRKFR IMEX scheme of~\cite{babbar2025crknoncons} to solve~\eqref{eq:lin.con.law.source}.
The $s$-stage scheme is specified by double \textit{Butcher tableaux}
\begin{equation}
  \begin{array}{c|c}
    \widetilde{\bc} & \tilde{A}\\
    \hline
    & \widetilde{\bb}^T
  \end{array}, \qquad \begin{array}{c|c}
    \bc & A\\
    \hline
    & \bb^T
  \end{array},
\label{eq:imex.butcher}
\end{equation}
where $\tilde{A} = (\tilde{a}_{ij})_{s \times s}$ is a strictly
lower triangular matrix corresponding to the explicit part.
The matrix $A =
(a_{ij})_{s \times s}$ corresponds to the implicit part.
We use diagonally implicit RK (DIRK) methods for the implicit part, i.e., the $A$ matrix is lower triangular.
The $\widetilde{\bc} = (\tilde{c}_1,
\ldots, \tilde{c}_s)^T, \bc = (c_1, \ldots, c_s)^T$~\cite{Hairer1991, Hairer1996} are the
\textit{quadrature nodes} given by $\tilde{c}_i = \sum_{j = 1}^{i - 1} \tilde{a}_{ij}, c_i = \sum_{j = 1}^i a_{ij}.$
The $\widetilde{\bb} = (\tilde{b}_1, \tilde{b}_2, \ldots, \tilde{b}_s)^T, \bb = (b_1, b_2, \ldots, b_s)^T$ are the \textit{weights} used to perform the final evolution with evaluations from the $s$ stages.
A discussion on choices of the Butcher tableaux for the cRKFR method can be found in Appendix C of~\cite{babbar2025crknoncons}.
Following~\cite{chen2024,babbar2025crk}, we define a local differentiation operator $\dlocx$ where the inter-element terms from~\eqref{eq:dfrx.defn} are dropped to be $\dlocx \pf (\uud_h) (\xi_p) = \pdx \discf (\xi_p)$.
The operator $\dlocx$ computes derivatives of the degree $N$ discontinuous flux approximation~\eqref{eq:discts.flux}.
In practice, this operation is performed with a differentiation matrix~\eqref{eq:diff.matrix}.
The idea of using local differentiation operators was also used in combination with Continuous Extension Runge-Kutta (CERK) methods~\cite{owren1991} to construct a high order local space-time predictor solution for ADER-DG (Arbitrary high order DERivative - Discontinuous Galerkin) schemes~\cite{Dumbser2008,Dumbser2014} in~\cite{gassner2011}.
The compact Runge-Kutta flux reconstruction method in the time average framework~\cite{babbar2025crk,babbar2022lax} applied to~\eqref{eq:lin.con.law.source} is given by
\begin{eqnarray}
\uu^{(i)} & = & \uu^n - \Delta t \sum_{j = 1}^{i - 1} \tilde{a}_{i j} \dlocx  \pf_h ( \uu^{(j)} ) +
\Delta t \sum_{j = 1}^i a_{ij}  \bss ( t_n + c_j \Delta t, \uu^{(j)} ), \label{eq:crk.fr.inner} \\
&& \qquad \qquad \qquad \qquad \qquad \qquad \qquad \qquad i = 1, \ldots, s, \label{eq:crkfr.inner} \nonumber\\
\uu^{n + 1} & = & \uu^n
- \Delta t \dfrx \F
+ \Delta t \bS_h^{\delta}, \label{eq:crkfr}\\
\dfrx \F &=& \pdx \F_h  = \pdx \discF + (\Fnum_\eph - \discFrefe(1))\pdx g_R + (\Fnum_\emh - \discFrefe(0)) \pdx g_L, \label{eq:cts.F}
\end{eqnarray}
where the time averaged flux $\F$ is computed in the reference coordinates as~\cite{babbar2025crk}
\begin{equation}
\discFref = \sum_{p = 0}^N \F_{e, p} \ell_p (\xi), \qquad \F_{e, p} =
\sum_{i = 1}^s \tilde{b}_i  \pf (\uu_{e, p}^{(i)}). \label{eq:disc.avg.flux}
\end{equation}
The time average source term $\bS_h^{\delta}$ and the time average solution $\uU_h^\delta$ are computed similarly~\cite{babbar2025crknoncons}.
The idea of time averaging is also inherently present in ADER-DG schemes~\cite{Dumbser2008,Dumbser2014} where the local space-time predictor solution is used to compute the time average flux and source term.
The time averages are used in~\cite{Dumbser2008,Dumbser2014} to perform evolution to the next time level, similar to~\eqref{eq:crkfr}.
A formal justification for why~\eqref{eq:disc.avg.flux} gives an approximation to the time averaged flux can be found in Appendix~A of~{\cite{babbar2025crk}}.
The \textit{time averaged numerical flux} $\Fnum_\eph$ used to construct the continuous flux approximation~\eqref{eq:cts.F} is computed as
\begin{equation} \label{eq:numflux.cons}
\Fnum_\eph = \half (\F_\eph^- + \F_\eph^+) - \frac{\sigma(\fmatrix)}{2} (\uU_\eph^+ - \uU_\eph^-),
\end{equation}
where $\F_\eph^\pm, \uU_\eph^\pm$ are the approximations of the time average flux and solution~(\ref{eq:disc.avg.flux}) at the interface $x_\eph$.
The dissipation term proportional to $\sigma(\fmatrix)$ is chosen such that the scheme has the same CFL number as the cRK scheme of~\cite{chen2024} while requiring only a single numerical flux $\Fnum_\eph$ to be computed per time step, see Section 3.1 of~\cite{babbar2025crk} for details.
The idea of using the time average solution to compute the dissipation term in the numerical flux was first introduced for Lax-Wendroff schemes in~\cite{babbar2022lax} where it was shown to improve the CFL numbers.
For the linear case, the scheme of~\cite{babbar2022lax} is equivalent to the ADER-DG scheme of~\cite{Dumbser2008,Dumbser2014} and thus has the same CFL numbers as the ADER-DG scheme~\cite{gassner2011}.
The difference between the cRKFR scheme~\eqref{eq:crkfr} and the standard RKFR scheme lies in the fact that the inner stages that evaluate $\{\uu_h^{(i)} \}_{i = 1}^s$ in~\eqref{eq:crkfr.inner} use the local operator $\dlocx$ instead of $\dfrx$ in~\eqref{eq:semidiscretization.dg}.

As was shown through various numerical experiments in~\cite{babbar2025crknoncons}, the cRKFR IMEX method described above is able to treat stiff source terms while requiring only the time step restrictions of the advective part of the system.

\section{Review of the Jin-Xin relaxation system} \label{sec:jin.xin}
The Jin-Xin relaxation system~\cite{jinxin1995} approximates a non-linear system of hyperbolic conservation laws by a system of hyperbolic balance laws with linear fluxes and stiff source terms.
Since there are some subtleties in the asymptotic analysis of the 2-D case, we present the analysis for the 2-D system of conservation laws
\begin{equation}
\partial_t \uu + \partial_{x_1} \pf_1 + \partial_{x_2} \pf_2 = \bzero. \label{eq:cons.law}
\end{equation}
Corresponding to each dimension, Jin-Xin relaxation introduces two new variables $\bv_1, \bv_2$ and a relaxation parameter $\varepsilon > 0$, leading to the larger system
\begin{equation}
\begin{split}
\partial_t  \uu + \partial_{x_1}  \bv_1 + \partial_{x_2}  \bv_2 & =
\boldsymbol{0},\\
\partial_t  \bv_1 + \ba_1^2 \partial_{x_1}  \uu & = -
\frac{1}{\varepsilon}  (\bv_1 - \pf_1 (\uu)) ,\\
\partial_t  \bv_2 + \ba_2^2 \partial_{x_2}  \uu & = -
\frac{1}{\varepsilon}  (\bv_2 - \pf_2 (\uu)),
\end{split} \label{eq:jin.xin}
\end{equation}
where $\ba_1 = a_1 \bI$, $\ba_2 = a_2 \bI$, $\bI$ is the identity matrix, and $a_1, a_2 > 0$ are constants that will be chosen later.
The second and third equations of~\eqref{eq:jin.xin} (formally) yield $\bv_i = \pf_i(\uu) + O(\varepsilon)$.
Inserting this into the first equation of~\eqref{eq:jin.xin}, we obtain
\begin{equation}
\partial_t \uu + \partial_{x_1} \pf_1 + \partial_{x_2} \pf_2 = \mathcal{O}(\varepsilon),
\label{eq:jin.xin.conv}
\end{equation}
so that~\eqref{eq:jin.xin} formally converges to~\eqref{eq:cons.law} as $\varepsilon \to 0$.
Thus, the Jin-Xin relaxation approximates~\eqref{eq:cons.law} by a system of balance laws with linear fluxes and stiff source terms, which can be solved numerically without a non-linear Riemann solver.

We now perform an asymptotic analysis to better understand the $\mathcal{O}(\varepsilon)$ term.
This analysis will help us understand the conditions on $a_1, a_2$ for the system~\eqref{eq:jin.xin} to be well-posed, and also how the relaxation parameter $\varepsilon$ relates to the numerical dissipation of the system.
This analysis has been performed in~\cite{jinxin1995,liu1987,driollet2004}.
The conditions for $a_i$ were obtained in~\cite{jinxin1995} in terms of the eigenvalues of the products of flux Jacobians and entropy Hessians.
To give an intuitive understanding, we present a simpler analysis yielding sufficient bounds in terms of the norms of the flux Jacobians.
The obtained bounds match the conditions obtained in~\cite{jinxin1995}, but the eigenvalues are replaced by the norms.
In the numerical results (Section~\ref{sec:num.results}), we choose $a_1, a_2$ based on eigenvalues of flux Jacobians, as they are easier to compute than the bounds based on the norms of Jacobians or on the eigenvalues of the products of flux Jacobians and entropy Hessians.
This asymptotic analysis is called the \emph{Chapman-Enskog expansion}~\cite{jinxin1995} of~\eqref{eq:jin.xin}, which is also used to formally derive the Navier-Stokes equation from the Boltzmann equation in the limit of small Knudsen number~\cite{chapman1990mathematical}.
The first step is to use the equations for $\bv_i$ in~\eqref{eq:jin.xin} to get the expansion
\begin{equation}
\bv_i = \pf_i(\uu) - \varepsilon ( \partial_t \bv_i + \ba_i^2 \partial_{x_i} \uu ) = \pf_i(\uu) + \mathcal{O}(\varepsilon).
\label{eq:v_expansion}
\end{equation}
Taking a temporal derivative of equation~\eqref{eq:v_expansion} gives us, using~\eqref{eq:jin.xin.conv},
\begin{align*}
\partial_t \bv_i &= \partial_t \pf_i(\uu) + \mathcal{O}(\varepsilon)
= \pf_i'(\uu) \partial_t \uu + \mathcal{O}(\varepsilon)
= - \pf_i'(\uu) (\partial_{x_1} \pf_1 + \partial_{x_2} \pf_2) + \mathcal{O}(\varepsilon)\\
&= -(\pf_i'(\uu) \pf_1'(\uu) \partial_{x_1} \uu + \pf_i'(\uu) \pf_2'(\uu) \partial_{x_2} \uu) + \mathcal{O}(\varepsilon).
\end{align*}
Substituting this result into the expansion of $\bv_i$~\eqref{eq:v_expansion} gives us
\begin{equation*}
\bv_i = \pf_i(\uu) - \varepsilon ( -(\pf_i'(\uu) \pf_1'(\uu) \partial_{x_1} \uu + \pf_i'(\uu) \pf_2'(\uu) \partial_{x_2} \uu) + \ba_i^2 \partial_{x_i} \uu ) + \mathcal{O}(\varepsilon^2).
\end{equation*}
Substituting the above expansion of $\bv_i$ in the first equation of~\eqref{eq:jin.xin}, we obtain
\begin{equation}
\begin{aligned}
& \partial_t \uu + \partial_{x_1} \pf_1 + \partial_{x_2} \pf_2 = \varepsilon \left[ \partial_{x_1} \left( (\ba_1^2 - \pf_1'(\uu)^2) \partial_{x_1} \uu - \pf_1'(\uu)\pf_2'(\uu) \partial_{x_2} \uu \right) \right. \\
& \qquad \qquad \qquad \qquad \qquad \left. + \partial_{x_2} \left( (\ba_2^2 - \pf_2'(\uu)^2) \partial_{x_2} \uu - \pf_2'(\uu)\pf_1'(\uu) \partial_{x_1} \uu \right) \right] + \mathcal{O}(\varepsilon^2).
\end{aligned}
\label{eq:jin.xin.diffusion}
\end{equation}
In order for the system to be well-posed, we would require the right-hand side to not be anti-diffusive. In the 1-D case, the condition is simply given by $a_1 \ge \|\pf_1'\|$.
The 2-D case requires us to study the second-order PDE operator
\begin{equation*}
\begin{gathered}
L = \partial_{x_1} \left( (\ba_1^2 - \fA_{11}) \partial_{x_1}- \fB_{12} \partial_{x_2} \right) + \partial_{x_2} \left( (\ba_2^2 - \fA_{22}) \partial_{x_2} - \fB_{21} \partial_{x_1} \right), \\
\fA_{ij} = \pf_i'(\uu) \pf'_j(\uu).
\end{gathered}
\end{equation*}
As is standard for such analysis, we assume periodic boundary conditions to be able to perform integration by parts without the appearance of boundary terms.
By multiplying the operator by a test function~$\boldsymbol{\phi}$, integrating over the physical domain and performing integration by parts, we obtain its bilinear form $B(\uu,\boldsymbol{\phi})$ given by
\[
B(\uu,\boldsymbol{\phi}) = - \int ( (\ba_1^2 - \fA_{11}) \partial_{x_1} \uu - \fB_{12} \partial_{x_2} \uu ) \partial_{x_1} \boldsymbol{\phi} + ( (\ba_2^2 - \fA_{22}) \partial_{x_2} \uu - \fB_{21} \partial_{x_1} \uu ) \partial_{x_2} \boldsymbol{\phi} \, \ud x.
\]
In order for the operator to be diffusive, we would require $B(\uu,\uu) \le 0$ for all $\uu$.
For convenience of notation, we define
\[
\|\fA_{ij}\| = \|\pf_i' \pf_j' \| := \sup_{\uu} \| \pf_i'(\uu) \pf_j'(\uu) \|_2,
\]
where $\| \cdot \|_2$ is the standard $L^2$ operator norm.
Assuming $a_i \ge  \| \pf_i'\|$, we have
\begin{multline*}
B(\uu, \uu) \le  -(a_1^2 - \| \fA_{11} \|) \| \partial_{x_1} \uu \|_2^2
- (a_2^2 - \| \fA_{22} \|) \| \partial_{x_2} \uu \|_2^2 \\
+ \| \fB_{12} \| \| \partial_{x_1} \uu \|_2 \| \partial_{x_2} \uu \|_2
+ \| \fB_{21} \| \| \partial_{x_1} \uu \|_2 \| \partial_{x_2} \uu \|_2.
\end{multline*}
Since $\| \fB_{i j} \| \le \| \pf_i' \| \| \pf_j' \|$,
the inequality will be satisfied if
\begin{equation*}
\mathcal{L} = \begin{bmatrix}a_1^2 - \|\pf_1'(\uu)\|^2 & - \|\pf_1'(\uu)\| \|\pf_2'(\uu)\| \\ -\|\pf_1'(\uu)\| \|\pf_2'(\uu)\| & a_2^2 - \|\pf_2'(\uu)\|^2 \end{bmatrix}
\end{equation*}
is positive semi-definite.
This requires the trace and determinant of the operator to be non-negative, which gives us the following conditions
\begin{equation}\label{eq:subcharacteristic}
\begin{split}
a_1^2 + a_2^2 & \ge \|\pf_1'(\uu)\|^2 +  \|\pf_2'(\uu)\|^2,\\
a_1^2 a_2^2  - a_1^2 \|\pf_2'(\uu)\|^2 - a_2^2 \|\pf_1'(\uu)\|^2 & \ge 0.
\end{split}
\end{equation}
The first condition is satisfied by choosing $a_i \ge \|\pf_i'\|$.
For the second condition, dividing by $a_1^2 a_2^2$ gives us the elliptic condition
\begin{equation}\label{eq:elliptical.condition}
\frac{\|\pf_1'(\uu)\|^2}{a_1^2} + \frac{\|\pf_2'(\uu)\|^2}{a_2^2} \le 1,
\end{equation}
which is the same condition as the one obtained in~\cite{jinxin1995}, but the norms are replaced by eigenvalues of the products of flux Jacobians and entropy Hessians.
The condition~\eqref{eq:elliptical.condition} is not satisfied just by choosing $a_i \ge \|\pf_i'\|$, but it requires choosing $a_i$ to be slightly larger.
For example, we could choose $a_i \ge \sqrt{2} \|\pf_i'\|$ or $a_i \ge \sqrt{\| \pf'_1 \|^2 + \|\pf'_2\|^2}$.
The former choice is used for the numerical experiments in Section~\ref{sec:num.results}.
In practice, we use eigenvalues of the flux Jacobians instead of norms, since computing norms is more expensive.
Although choosing $a_i$ based on eigenvalues does not guarantee that the system will be dissipative, we found that the choice works well in practice for the problems we considered.

\section{Shock capturing with Jin-Xin relaxation} \label{sec:shock.capturing}


Choosing $\ba_i$ in~\eqref{eq:jin.xin} appropriately, equation~\eqref{eq:jin.xin.diffusion} shows that the Jin-Xin relaxation system can be interpreted as adding a diffusion term proportional to $\varepsilon$ to the original system of conservation laws~\eqref{eq:cons.law}, ignoring $\mathcal{O}(\varepsilon^2)$ terms.
Thus, by choosing $\varepsilon$ to be large in the vicinity of shocks and discontinuities, and small in smooth regions, we can use the Jin-Xin relaxation system as a shock-capturing method for the original system of conservation laws.
Therefore, we will be able to use the high-order cRKFR method to solve hyperbolic conservation laws with shocks and discontinuities without the need for a non-linear Riemann solver.
Since the cRKFR IMEX method solves the implicit equation~\eqref{eq:crk.fr.inner} locally at each solution point, the scheme will be prescribed once we specify the value of $\varepsilon$ in each element.
In order to do that, we use the smoothness indicator of~\cite{hennemann2021} which is inspired from the idea of~\cite{Persson2006}.
The degree $N$ solution polynomial~\eqref{eq:soln.poly} can be written in terms of an orthogonal basis like Legendre polynomials.
The smoothness of the solution can be assessed by analyzing the decay of the coefficients of the orthogonal expansion, a technique originally proposed by Persson and Peraire~\cite{Persson2006} and subsequently refined by Henemann et al.~\cite{hennemann2021}.

Let $q = q(\uu)$ be the quantity used to measure the solution smoothness. We first project this onto Legendre polynomials $q_h(\xi) = \sum_{j=0}^N \hat q_j L_j(2\xi-1)$, for $\xi\in [0,1]$, where $\hat q_j = \int_0^1 q(\uu_h(\xi)) L_j(2\xi-1) \ud\xi$.
The Legendre coefficients $\hat q_j$ are computed using the collocation quadrature as $\hat q_j = \sum_{q=0}^N q(\uu_q^e) L_j(2\xi_q-1) w_q$.
Then the energy contained in the highest modes relative to the total energy of the polynomial is computed as $\en = \max \left(
\frac{\hat q_{N-1}^2}{\sum_{j=0}^{N-1}\hat q_j^2},
\frac{\hat q_N^2}{\sum_{j=0}^N \hat q_j^2}
\right)$.
The $N^\text{th}$ Legendre coefficient $\hat q_N$ of a function that is in the Sobolev space $H^p$ decays as $\mathcal{O}(1/N^p)$ (see Chapter 5, Section 5.4.2 of~\cite{Canuto2007}).
Thus, for a smooth function, the energy $\en$ will be small, while for a non-smooth function, the energy $\en$ will be large.
In~\cite{hennemann2021}, the $\en$ value is carefully mapped to $\alpha \in [0,1]$ to be used as a coefficient to \textit{blend} with a lower-order method.
In this work, we map the $\en$ value to $\varepsilon$ to be used as a coefficient for the Jin-Xin relaxation system directly as
\[
\tilde{\varepsilon} = k \en,
\]
where $k$ is a parameter in our method; we always choose $k = 2 \times 10^5$.
Additionally, we set threshold values $\varepsilon_{\min}, \varepsilon_{\max}$ and clip the value of $\tilde{\varepsilon}$ to finally choose $\varepsilon$ as
\begin{equation}
\varepsilon = \min\{\varepsilon_{\max}, \max\{\tilde{\varepsilon}, \varepsilon_{\min}\}\}.
\label{eq:eps.thresholds}
\end{equation}
For the numerical tests, we needed to do multiple experiments to fine-tune the value of $\varepsilon_{\max}$ for every problem to get the best results, while we were able to make the natural choice of $\varepsilon_{\min} = 10^{-12}$ for all problems.
A very small value of $\varepsilon_{\max}$ lead to insufficient dissipation and spurious oscillations near shocks, while a very large value of $\varepsilon_{\max}$ lead to excessive dissipation and smearing of the solution near shocks.
This problem of fine-tuning $\varepsilon_{\max}$ is thus a limitation of our method that needs to be addressed in future work.
It could be addressed by using a better mapping from $\en$ to $\varepsilon$, or a better smoothness indicator to compute $\en$.

\section{Numerical results} \label{sec:num.results}
A comparison of the Jin-Xin relaxation based shock-capturing method with the subcell-based blending limiter of~\cite{babbar2024admissibility,babbar2025crk,hennemann2021} is performed for test problems involving scalar equations and the compressible Euler equations.
The schemes presented use Gauss-Legendre (GL) solution points with Radau correction functions~\eqref{eq:cts.F}~\cite{Huynh2007}.
This choice is known to show more accuracy in comparison to other choices of solution points and correction functions~\cite{Huynh2007,babbar2022lax,babbar2024admissibility}, while alternative choices can improve stability, e.g., by increasing the CFL stability limit.
The proposed scheme is general, and can be used with other choices of solution points and correction functions as well.
Unless specified otherwise, the subcell-based blending limiter performs first-order blending on the subcells, although higher-order blending can also be performed~\cite{babbar2024admissibility,babbar2025crk,ramirez2022}.

The code~{\tt Tenkai.jl}~\cite{tenkai} is used to perform all the numerical experiments in this work, and the code to reproduce the results is available at~\cite{babbar2026jinxinrepro}.
\subsection{Scalar equations}
\subsubsection{Burgers' equation: convex flux}
We test the Jin-Xin relaxation system for the Burgers' flux function $f(u) = \frac{u^2}{2}$, with the initial condition $u(x,0) = 2 + \sin(\pi (x - 0.7))$ for $x \in [-1, 1]$ with periodic boundary conditions with final time $t = 0.5$.
A value of $\varepsilon_{\max} = 2 \times 10^{-3}$ is used for the Jin-Xin shock-capturing method, with SSP3-IMEX(4,3,3) scheme of~\cite{pareschi2005} for the time integration.
The results are shown in Figure~\ref{fig:scalar}(a) on a grid of 20 elements with polynomial degree $N = 3$.
The subcell-based blending limiter and Jin-Xin relaxation based shock-capturing method give similar results for this problem, with the Jin-Xin method being slightly less dissipative than the subcell-based blending limiter.
\subsubsection{Buckley-Leverett equation: non-convex flux}
\begin{figure}
\centering
\begin{tabular}{cc}
\includegraphics[width=0.46\columnwidth]{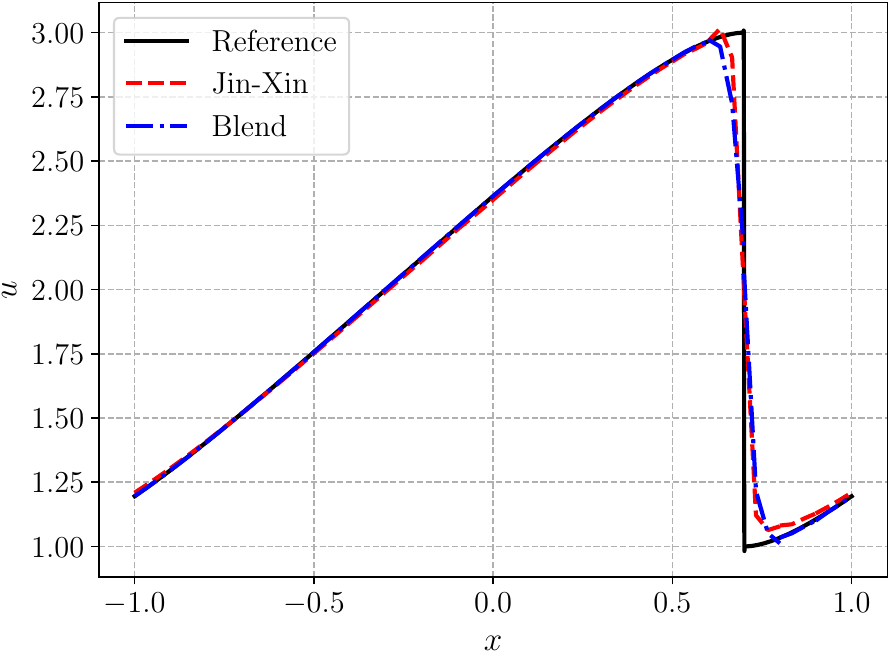} &
\includegraphics[width=0.46\columnwidth]{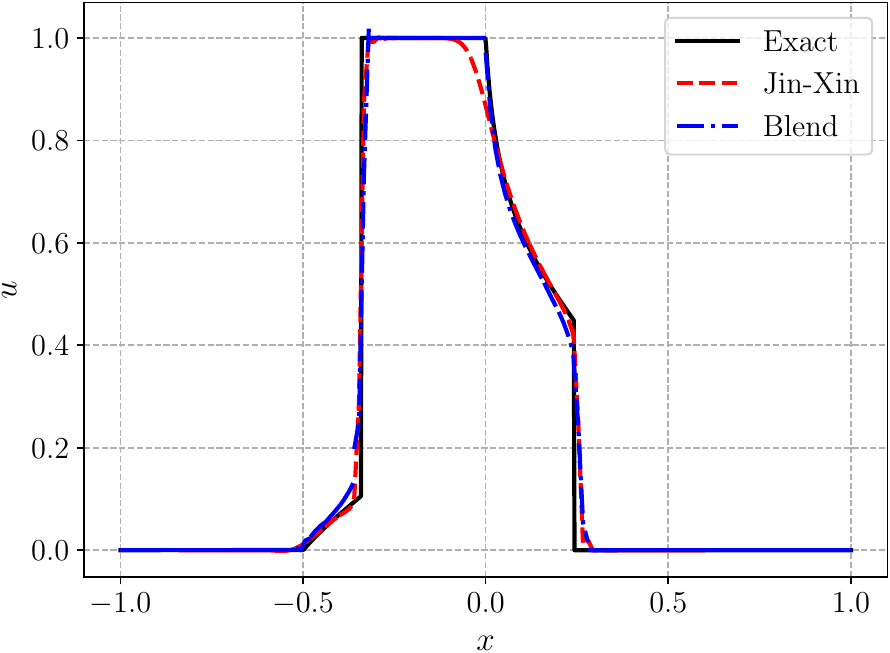} \\
(a) & (b)
\end{tabular}
\caption{Comparison of subcell-based blending limiter and the Jin-Xin relaxation based limiter for (a) Burgers' equation and (b) Buckley-Leverett equation. The SSP3-IMEX(4,3,3) scheme of~\cite{pareschi2005} is used for the Jin-Xin relaxation method.}
\label{fig:scalar}
\end{figure}

We consider the Buckley-Leverett equation $u_t+f(u)_x=0,$ where the flux $f(u)=\dfrac{4u^2}{4u^2+(1-u)^2}$ is convex and concave in different regions of the solution space. The numerical solution is computed up to the time $t=0.15$ with the initial condition $u(x,0) = \chi_{[-\frac 12, 0]}$, where $\chi$ is the indicator function.
The physical domain is $[-1,1]$ and periodic boundary conditions are used.
The threshold value $\varepsilon_{\max}$ is set to $2 \times 10^{-4}$ for the Jin-Xin shock-capturing method, and the SSP3-IMEX(4,3,3) scheme of~\cite{pareschi2005} is used for the time integration.
The results are shown in Figure~\ref{fig:scalar}(b) on a grid of 50 elements with polynomial degree $N = 3$.
The Jin-Xin relaxation based shock-capturing method is able to capture the shock and rarefaction waves in the solution, although it is more dissipative than the subcell-based blending limiter for this problem.

\subsection{Euler's equations}
\subsubsection{Woodward-Colella blast wave problem}
\begin{figure}
\centering
\begin{tabular}{cc}
\includegraphics[width=0.46\columnwidth]{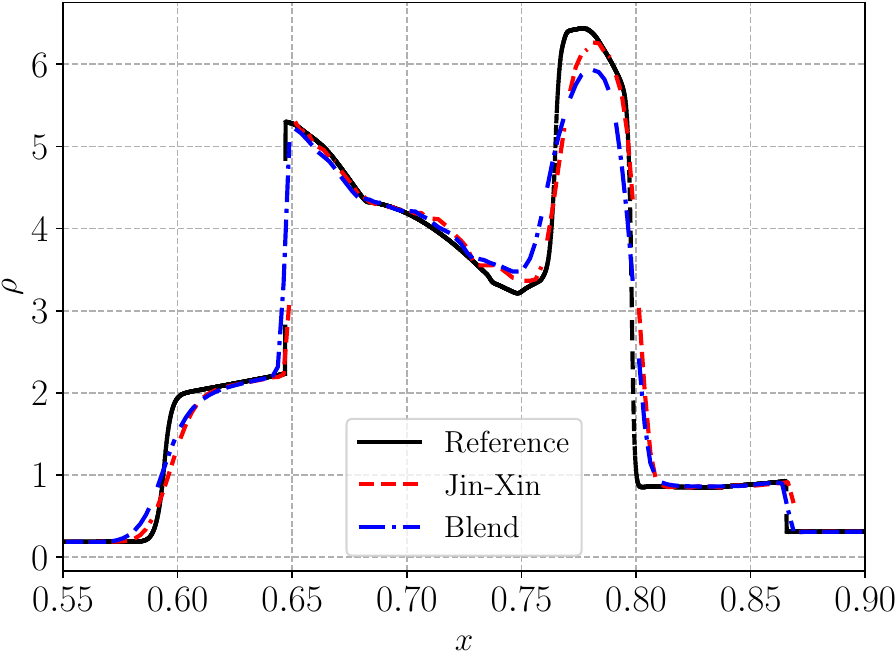} &
\includegraphics[width=0.46\columnwidth]{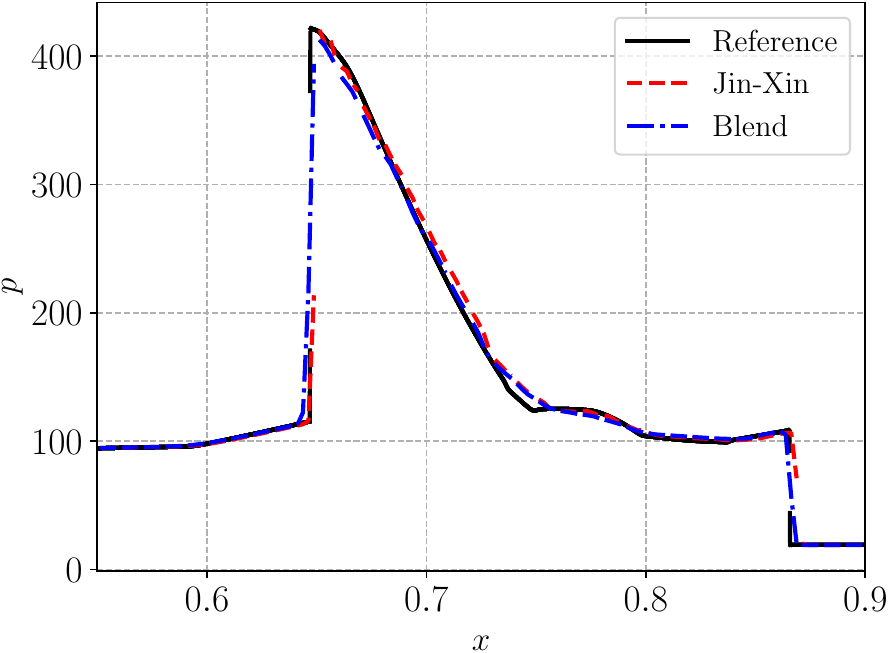} \\
(a) & (b)
\end{tabular}
\caption{Comparison of subcell-based blending limiter and the Jin-Xin relaxation based limiter for the blast wave problem showing (a) Density and (b) Pressure. The SSP3-IMEX(4,3,3) scheme of~\cite{pareschi2005} is used for the Jin-Xin relaxation method.}
\label{fig:blast}
\end{figure}

This blast wave test from~\cite{Woodward1984} consists of an initial condition with two discontinuities given as
\begin{equation*}
(\rho,v_1,p)=\begin{cases}
(1,0, 1000), & \mbox{ if } x<0.1,\\
(1,0,0.01), & \mbox{ if }  0.1 < x <0.9,\\
(1,0, 100), & \mbox{ if } x> 0.9,
\end{cases}
\end{equation*}
in the domain $[0,1]$. The boundaries are set as solid walls by imposing the reflecting boundary conditions at $x=0$ and $x=1$.
The two discontinuities lead to two Riemann problems, each of which leads to a shock, rarefaction and a contact discontinuity.
The boundary conditions cause reflection of shocks and expansion waves off the solid wall and several wave interactions inside the domain.
As time goes on, the two reflected Riemann problems interact, which is the main point of interest for this test~\cite{Woodward1984}.
The positivity limiter of~\cite{zhang2010c} is needed for both the subcell-based blending limiter and the Jin-Xin relaxation based shock-capturing method to ensure positivity of density and pressure.
The threshold value $\varepsilon_{\max}$ is set to $10^{-7}$ for the Jin-Xin shock-capturing method, and the SSP3-IMEX(4,3,3) scheme of~\cite{pareschi2005} is used for the time integration.
The results are shown in Figure~\ref{fig:blast} on a grid of 400 elements with polynomial degree $N = 3$.
The solutions look similar for both methods, although the Jin-Xin relaxation based shock-capturing method is slightly more accurate than the subcell-based blending limiter for capturing the smooth extremum.

\subsubsection{Sedov's blast problem}

\begin{figure}
\centering
\includegraphics[width=0.7\columnwidth]{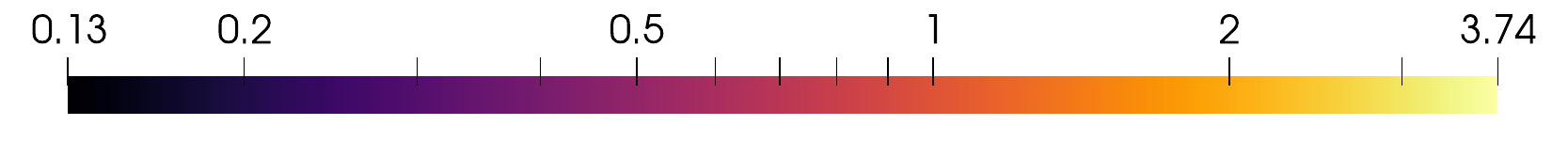} \\
\begin{tabular}{cc}
\includegraphics[width=0.44\columnwidth]{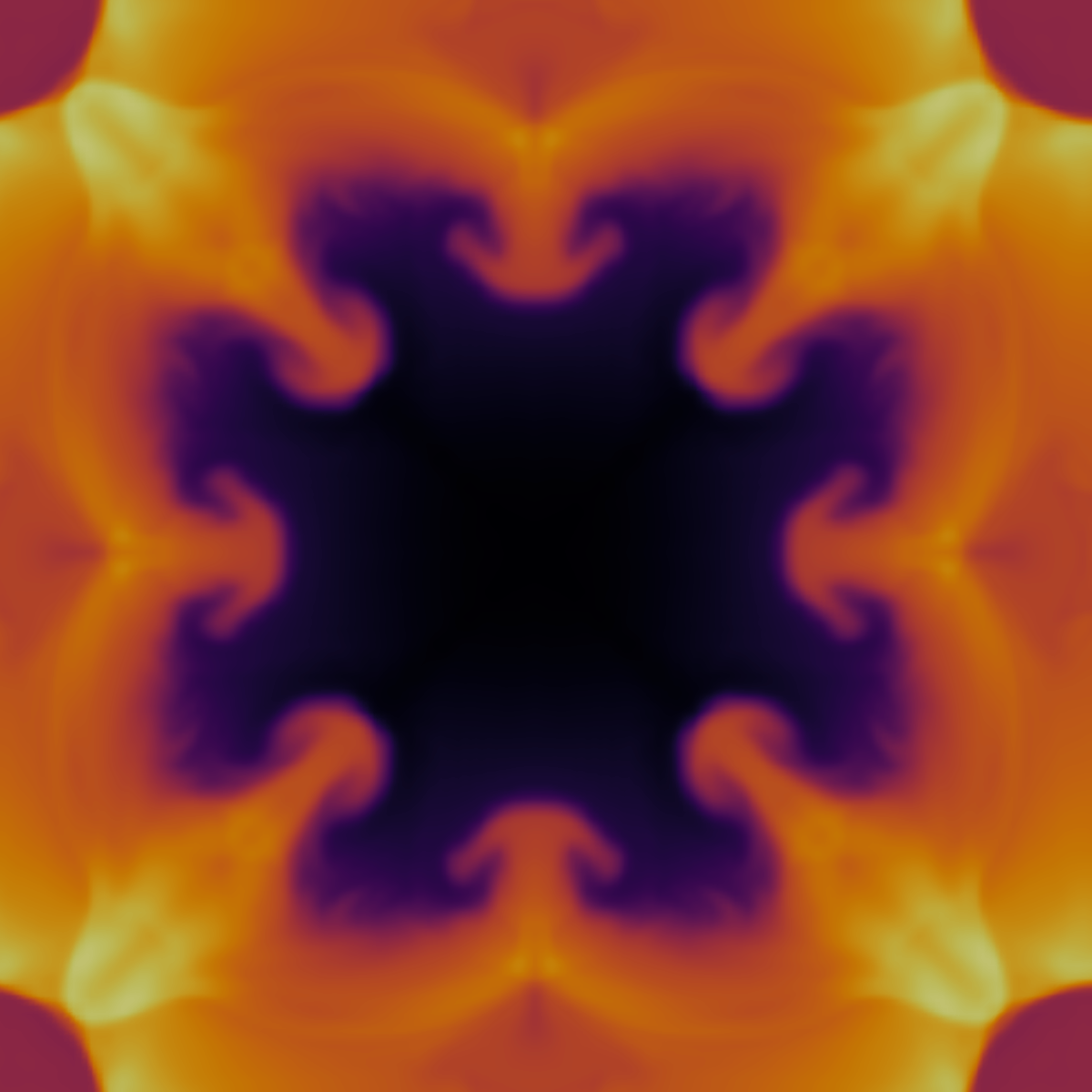} &
\includegraphics[width=0.44\columnwidth]{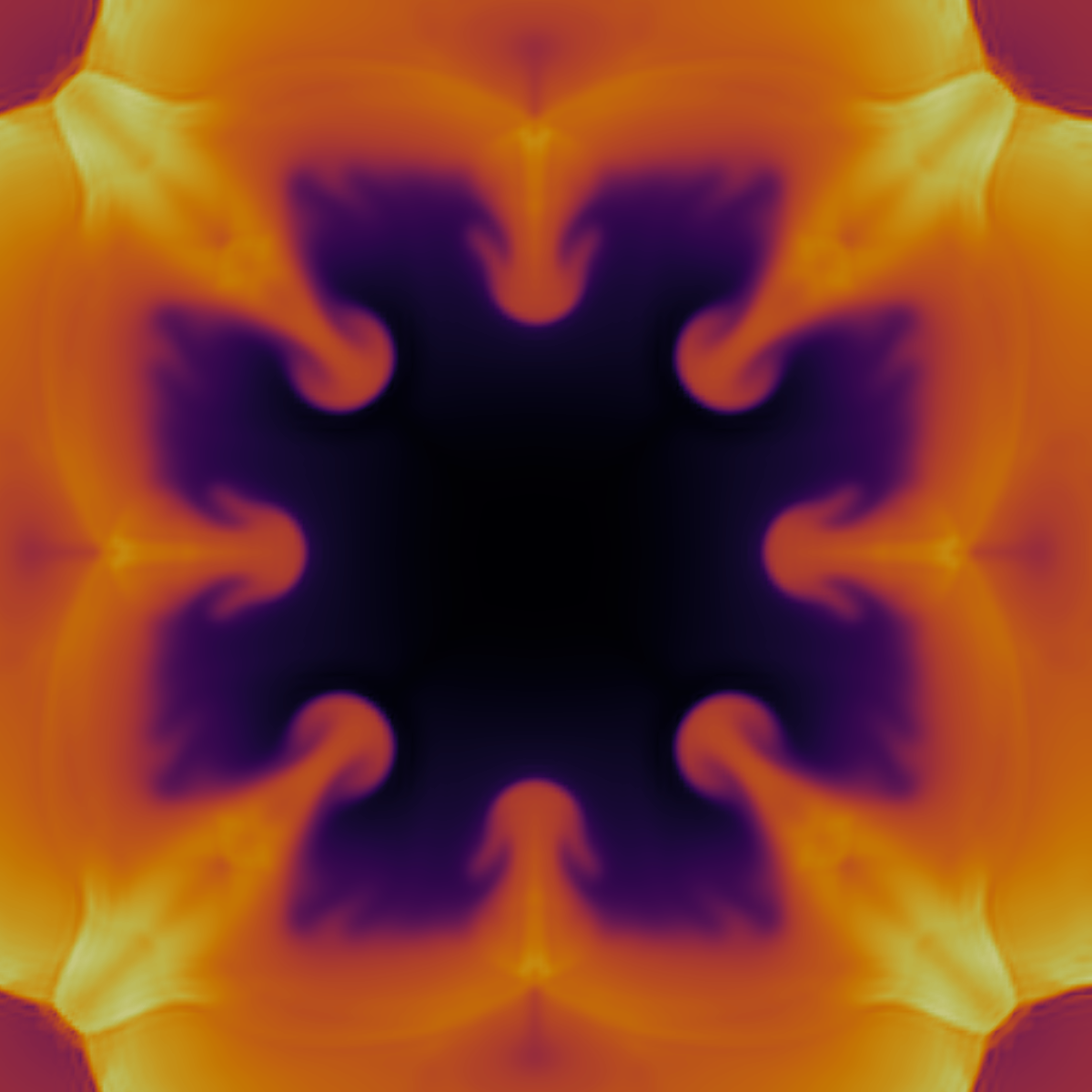} \\
(a) subcell-based blending scheme & (b) Jin-Xin relaxation based shock capturing scheme
\end{tabular}
\caption{Comparison of subcell-based blending limiter and the Jin-Xin relaxation based shock capturing scheme for the Sedov's blast wave problem. The BPR(3, 4, 3) scheme of~\cite{boscarino2013} is used for the Jin-Xin relaxation method.}
\label{fig:sedov.blast}
\end{figure}

This Sedov's blast test was introduced in~\cite{ramirez2021} to demonstrate the positivity-preserving property of their scheme.
This test consists of energy concentrated at the origin.
More precisely, for the initial condition, we assume that the gas is at rest ($v_1 = v_2 = 0$) with Gaussian distribution of density and pressure
\begin{equation}
\rho(x,y) = \rho_0 + \frac{1}{4\pi\sigma_\rho^2} \exp \left( -\frac{r^2}{2\sigma_\rho^2} \right), \ p(x,y) = p_0  + \frac{\gamma - 1}{4 \pi \sigma_p^2} \exp\left( -\frac{r^2}{2\sigma_p^2} \right), \ r^2 = x^2 + y^2,
\end{equation}
where $\sigma_\rho = 0.25$ and $\sigma_p = 0.15$.
The ambient density and ambient pressure are set to $\rho_0 = 1$, $p_0 = 10^{-5}$.
The boundaries are taken to be periodic, leading to interaction between shocks as they re-enter the domain through periodic boundaries and form small scale fractal structures.
The threshold value $\varepsilon_{\max}$ is set to $8 \times 10^{-5}$ for the Jin-Xin shock-capturing method, and the BPR(3, 4, 3) scheme of~\cite{boscarino2013} is used for the time integration.
The results are shown in Figure~\ref{fig:sedov.blast} on a grid of $64 \times 64$ elements with polynomial degree $N = 3$.
In this test, the subcell-based blending limiter uses MUSCL-Hancock reconstruction on the subcells~\cite{babbar2024admissibility,babbar2025crk}, and the results are very similar for both methods.
The differences are very slight, but the Jin-Xin relaxation based shock-capturing method is closer to the reference solutions shown in~\cite{babbar2024admissibility,babbar2025crk}.

\subsubsection{Kelvin-Helmholtz instability}
\begin{figure}
\centering
\includegraphics[width=0.7\columnwidth]{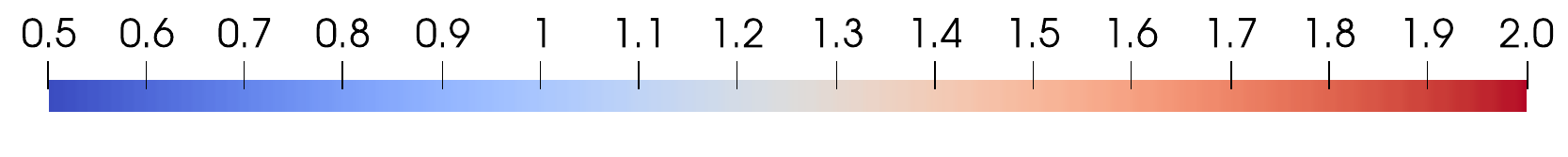} \\
\begin{tabular}{cc}
\includegraphics[width=0.44\columnwidth]{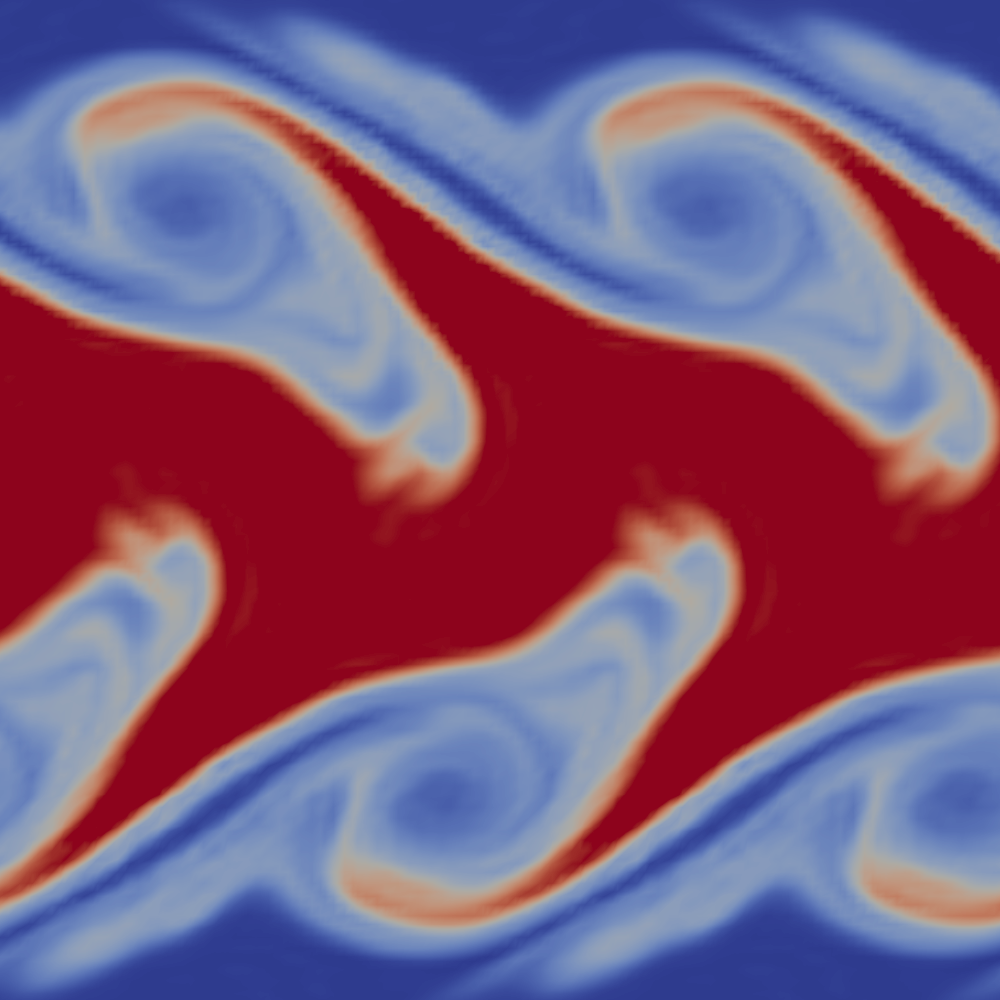} &
\includegraphics[width=0.44\columnwidth]{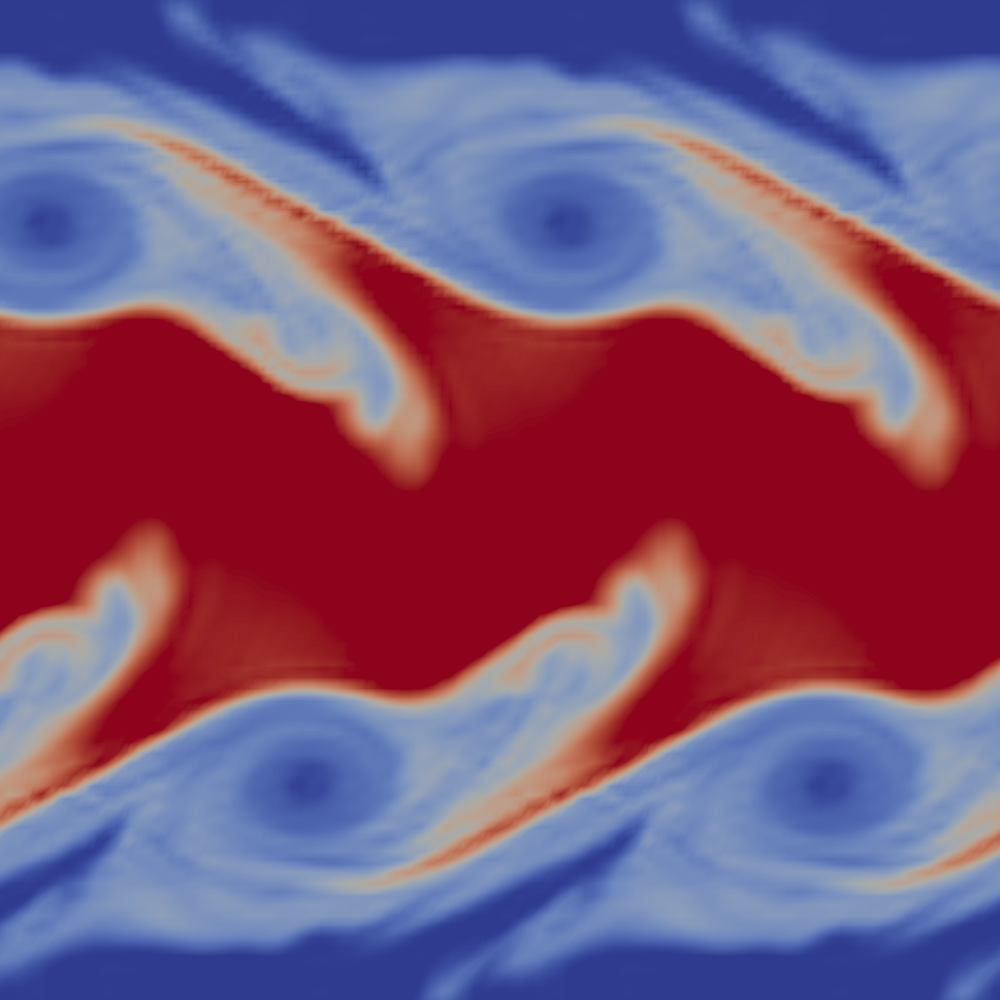} \\
(a) subcell-based blending scheme & (b) Jin-Xin relaxation based shock capturing scheme
\end{tabular}
\caption{Comparison of subcell-based blending limiter and the Jin-Xin relaxation based shock capturing scheme for the Kelvin-Helmholtz instability problem. The BPR(3, 4, 3) scheme of~\cite{boscarino2013} is used for the Jin-Xin relaxation method.}
\label{fig:khi}
\end{figure}

Fluid instabilities are essential for mixing processes and turbulence production, and play a significant role in many phenomena such as astrophysics.
The Kelvin-Helmholtz instability is a common fluid instability that occurs across contact discontinuities in the presence of a tangential shear flow.
We take the setup of the Kelvin-Helmholtz instability problem from~\cite{ramirez2021} which is given by $(\rho, p, u, v) = (\frac 12 + \frac 34 B, 1, \frac 12 (B-1), \frac{1}{10} \sin(2 \pi x))$,
where $B(x,y)= \tanh(15y + 7.5) - \tanh(15y - 7.5)$ is a smoothed approximation to a discontinuous step function on the domain $[-1,1] \times [-1,1]$ with periodic boundary conditions.
The simulation is run until time $t=10$ when the solution has developed pronounced vortex structures.
The threshold value $\varepsilon_{\max}$ is set to $10^{-6}$ for the Jin-Xin shock-capturing method, and the BPR(3, 4, 3) scheme of~\cite{boscarino2013} is used for the time integration.
The results are shown in Figure~\ref{fig:khi} on a grid of $32 \times 32$ elements with polynomial degree $N = 3$.
This test does not have strong discontinuities, and thus an important question about shock capturing schemes is how much additional dissipation is needed to run the simulation without requiring the positivity limiter~\cite{zhang2010c}.
For the subcell-based blending limiter, we could run the simulation by limiting the blending coefficient to be 0.01.
The Jin-Xin relaxation shock capturing scheme captures more small scale structures than the subcell-based blending limiter, which gives us an indication that it requires less additional dissipation to run the simulation.
Although there are some oscillations and additional artifacts, we have observed that they can be eliminated by slightly increasing the value of $\varepsilon_{\max}$, but we chose to show results with this value to demonstrate that the method can run with dissipation that is apparently lesser than the subcell-based blending limiter.
\section{Conclusions and future work} \label{sec:conclusion}
We introduced a Jin-Xin relaxation based shock-capturing method for high-order DG/FR methods for hyperbolic conservation laws.
The stiff source term of the Jin-Xin relaxation system is solved using the compact Runge-Kutta flux reconstruction (cRKFR) IMEX method of~\cite{babbar2025crknoncons}, which is an arbitrarily high-order method for solving hyperbolic equations with stiff source terms.
Numerical experiments are performed for scalar equations and the compressible Euler equations involving discontinuities and small-scale features.
The results show that the method is of comparable accuracy and robustness to the subcell-based blending limiter of~\cite{hennemann2021,babbar2024admissibility,babbar2025crk}.
Thus, a high-order robust shock-capturing method is obtained without the need for a non-linear Riemann solver.
Given its simplicity and comparable accuracy to pre-existing methods, further study of the method including a performance study and extension to more complex systems will be important research directions.
The issue with the method is that it currently requires fine-tuning of the threshold parameter $\varepsilon_{\max}$ for each problem; overcoming this will be an interesting future research direction.
This work was mostly limited to periodic problems, and proper treatment of general boundary conditions for the Jin-Xin relaxation model will also be important.
The last, but rather difficult open question is about provable fully discrete asymptotic preservation of the method, as $\varepsilon$ tends to zero.

\section*{Acknowledgments}

MA, AB and HR were supported by the Deutsche Forschungsgemeinschaft
(DFG, German Research Foundation, project number 528753982
as well as within the DFG priority program SPP~2410 with project number 526031774).
AB was also supported by the Alexander von Humboldt Foundation.
MSL acknowledges funding by the DFG - project number 463312734; 528753982.
GG acknowledges funding by the DFG under Germany's Excellence Strategy - EXC 3037- 533607693. GG further acknowledges funding through the German Federal Ministry for Education and Research (BMFTR) project “ICON-DG” (01LK2315B) of the “WarmWorld Smarter” program and received funding through the DFG research unit FOR 5409 “SNUBIC”.

\bibliographystyle{spbasic}
\bibliography{references.bib}
\end{document}